\documentclass[onecolumn]{autart}    

\usepackage{graphicx}          
\usepackage{amsmath,amssymb}
\usepackage{tikz}
\usepackage{comment}
\usepackage[round]{natbib}
\usetikzlibrary{shapes,arrows}
\begin{document}

\begin{frontmatter}

\title{A constraint-separation principle in model predictive control\thanksref{footnoteinfo}} 

\thanks[footnoteinfo]{This paper was not presented at any IFAC 
meeting. Corresponding author U.~V.~Kalabi\'{c}.}

\author[merl]{Uro\v{s} V.~Kalabi\'{c}}\ead{kalabic@umich.edu}~
\author[umich]{Ilya V.~Kolmanovsky}\ead{ilya@umich.edu}               

\address[merl]{Mitsubishi Electric Research Laboratories, 201 Broadway, Cambridge, MA 02139}  
\address[umich]{Department of Aerospace Engineering, University of Michigan, 1320 Beal Avenue, Ann Arbor, MI 48109}             


\begin{abstract}                          
In this brief, we consider the constrained optimization problem underpinning model predictive control (MPC). We show that this problem can be decomposed into an unconstrained optimization problem with the same cost function as the original problem and a constrained optimization problem with a modified cost function and dynamics that have been precompensated according to the solution of the unconstrained problem. In the case of linear systems subject to a quadratic cost, the unconstrained problem has the familiar LQR solution and the constrained problem reduces to a minimum-norm projection. This implies that solving linear MPC problems is equivalent to precompensating a system using LQR and applying MPC to penalize only the control input. We propose to call this a constraint-separation principle and discuss the utility of both constraint separation and general decomposition in the design of MPC schemes and the development of numerical solvers for MPC problems.
\end{abstract}

\end{frontmatter}

\section{Introduction}

Model predictive control (MPC) is an optimization-based framework for determining constraint-admissible, stabilizing control inputs to open-loop control systems \citep{rawlings_book}. Conventionally, MPC is applied according to the schematic in Fig.~\ref{fig:MPC_ol} and presented as a constraint-enforcing, feedback control scheme, 
which can simultaneously stabilize a system and enforce constraints on that system.
This is in contrast to other constraint-enforcing schemes, such as reference governors \citep{rg_survey}, which are only used to enforce constraints in precompensated systems, and anti-windup schemes \citep{anti_windup} and barrier-function methods \citep{barrier_fn}, which are used to modify stabilizing control designs in order to enforce constraints.

In this paper, we propound the perspective that 
MPC is not substantially different from the alternative schemes. This is because the optimization problem solved by MPC can be decomposed into two separate optimization problems, the solution to one ensuring stability of the inner-loop, and the solution to the other computing an outer-loop modification that enforces constraints.
In the case of linear systems subject to a quadratic penalty, the inner-loop compensation has the familiar closed-form solution of the discrete-time linear-quadratic regulator (LQR) \citep{dorato71} and the constraint-enforcing, outer-loop problem is reduced to a minimum-norm projection, 
without a terminal cost term.

The above implies that linear MPC can be interpreted as an add-on, constraint enforcing mechanism similar to the extended command governor (ECG) \citep{gilbertong11}.
This was first shown by \citet{kouvaritakis2000} for the time-invariant case, i.e., where the dynamics are time-invariant and the terminal cost is obtained as the solution to the discrete-time, infinite-horizon LQR problem, and the benefits of this approach for numerical implementation were discussed by \citet{kouvaritakis02,rossiter_book,rossiter10}. Here we show the same result in the case of linear time-varying systems, which is important 
because numerical strategies for nonlinear MPC problems often involve sequentially solving linear time-varying MPC problems. 
We refer to the result as a constraint-separation principle and discuss its implications.

In addition to the above result, we consider decomposition in the general nonlinear setting. We show that constraint separation is not generally possible in the case of nonlinear MPC problems, since the initial MPC problem is decomposed into an unconstrained problem requiring a convenient closed-form solution, which does not generally exist, and another MPC problem that does not necessarily have a convenient structure. Nevertheless, we show that if the MPC problem has a locally linear-quadratic structure, then constraint separation holds locally.

The paper is structured as follows. Section 2 derives the decomposition of nonlinear MPC into stabilizing and constraint-enforcing optimization problems. Section 3 derives the constraint-separation principle for locally linear MPC problems. Section 4 discusses the implications of the constraint-separation principle. Section 5 is the conclusion.


\tikzstyle{block} = [draw, fill=white, rectangle, 
    minimum height=3em, minimum width=6em, font=\sffamily]
\tikzstyle{sum} = [draw, fill=white, circle]
\tikzstyle{input} = [coordinate]
\tikzstyle{output} = [coordinate]
\tikzstyle{pinstyle} = [pin edge={to-,thin,black}]

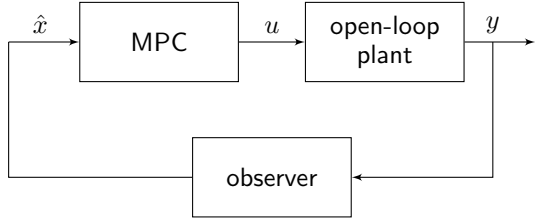
\begin{figure}[tbp]
\centering
\begin{tikzpicture}[auto, node distance=2cm,>=latex']
    \node [input, name=input] {};
    \node [block, right of=input] (controller) {MPC};
    \node [block, right of=controller,
            node distance=3cm,text width=0.7in,align=center] (system) {open-loop plant};
    \draw [->] (controller) -- node[name=u] {$u$} (system);
    \node [output, right of=system] (output) {};
    \node [block, below of=u] (measurements) {observer};
    
    \draw [->] (input) -- node [pos=0.45] {$\hat x$} (controller);
    \draw [->] (system) -- node [name=y,pos=0.4] {$y$}(output);
    \draw [->] (y) |- (measurements);
    \draw [-] (measurements) -|  node [near end] {} (input);
\end{tikzpicture}
\caption{MPC applied to an open-loop plant}\label{fig:MPC_ol}
\end{figure}

\section{MPC optimization problem}

The optimization problem we consider is given by \cite{rawlings_book},
\begin{subequations}\label{equ:mpc_prob}
\begin{align}
\min_u &~ V_f(x_N) + \sum_{k = 0}^{N-1} L_k(x_k,u_k), \label{equ:cost} \\
\text{sub.~to} &~x_{k+1} = f_k(x_k,u_k), \\
&~(x_k,u_k) \in \mathcal C_k,~\forall k \in \mathbb Z_N, \label{equ:ineq1}  \\
&~x_N \in \mathcal X_N, \label{equ:ineq2}
\end{align}
\end{subequations}
where $x_0$ is given, $f_k: \mathbb R^n \times \mathbb R^m \to \mathbb R^n$ is continuous, $f_k(0,0) = 0$, and $\mathcal C_k$ is closed for all $k\in\mathbb Z_N$ where $\mathbb Z_N$ is the set of the first $N$ non-negative integers. The cost functions $V_f:\mathbb R^n \to \mathbb R$ and $L_k:\mathbb R^n \times \mathbb R^m \to \mathbb R$ are continuous and locally bounded, and satisfy the following properties,
\begin{align*}
&L_k(0,0) = V_f(0) = 0, \\
&L_k(x,u) \geq \alpha(\|u\|),~V_f(x) \geq 0,
\end{align*}
for all $x \in \mathbb R^n$, $u \in \mathbb R^m$, $k\in \mathbb Z_N$, where $\alpha$ is a $\mathcal K_\infty$-function \citep{rawlings_book}. The assumptions are required 
to ensure the existence of a solution to the optimization problem. 
We note that our subsequent results do not dependent on other, standard assumptions found in the MPC literature \citep{rawlings_book}, like those which ensure that the solution is recursively feasible and results in a stabilizing control law.
We also note that we make no additional assumptions on the geometry of the sets $\mathcal C_k$ such as, for example, convexity. 

We define the sets,
\begin{align*}
&\bar{\mathcal C}_k := \{(x,u): (x,u) \in \mathcal C_k,~f_k(x,u) \in \mathcal X_{k+1}\}, \\
&\bar{\mathcal C}_k(x) := \{u: (x,u) \in \bar{\mathcal C}_k\},
\end{align*}
for all $k \in \mathbb Z_N$, where $\mathcal X_k := \operatorname{Proj}_{\mathbb R^n}{\mathcal C}_k$.
These sets are closed for all $x \in \mathbb R^n$ due to the closedness of $\mathcal C_k$ and $\mathcal X_{k+1}$ and the continuity of $f_k$.

A sequence of control inputs $u^*$ 
solving \eqref{equ:mpc_prob} satisfies \citep{rawlings_book},
\begin{subequations}
\begin{equation}\label{equ:u_opt}
u_k^* \in \underset{u \in \bar{\mathcal C}_k(x_k^*)}{\arg\min}~ V_{k+1}(f_k(x_k^*,u)) + L_k(x_k^*,u), \\
\end{equation}
for all $k \in \mathbb Z_N$, with $x_{k+1}^* = f_k(x_k^*,u_k^*)$  and $x_0^* = x_0$, where $V_k$ satisfies the Bellman equation,
\begin{equation}\label{equ:bell1}
V_k(x) = \min_{u \in \bar{\mathcal C}_k(x)} V_{k+1}(f_k(x,u)) + L_k(x,u),
\end{equation}
\end{subequations}
with domain $\bar{\mathcal X}_k := \operatorname{Proj}_{\mathbb R^n}\bar{\mathcal C}_k$, and $V_N = V_f|_{\mathcal X_N}$.

Consider the optimization problem \eqref{equ:mpc_prob} without inequality constraints \eqref{equ:ineq1}-\eqref{equ:ineq2}. A sequence of control inputs $\tilde u^*$ 
solving this problem satisfies, 
\begin{subequations}
\begin{equation}\label{equ:u_opt_unc}
\tilde u_k^* \in \underset{u}{\arg\min}~ \tilde V_{k+1}(f_k(\tilde x_k^*,u)) + L_k(\tilde x_k^*,u), \\
\end{equation}
for all $k \in \mathbb Z_N$, with $\tilde x_{k+1}^* = f_k(\tilde x_k^*,\tilde u_k^*)$ and $\tilde x_0^* = x_0$, where $\tilde V_k$ satisfies the Bellman equation,
\begin{equation}
\tilde V_k(x) = \min_{u}~\tilde V_{k+1}(f_k(x,u)) + L_k(x,u), 
\end{equation}
\end{subequations}
and $\tilde V_N = V_f$.
We call $\tilde V_k$ an unconstrained value function and distinguish it from the corresponding value function $V_k$.
A solution to the unconstrained problem defines a feedback law $\kappa_k:\mathbb R^n \to \mathbb R^m$ satisfying,
\begin{equation}
\kappa_k(\tilde x_k^*) = \tilde u_k^*, 
\end{equation}
where $\tilde u_k^*$ is a minimizer of \eqref{equ:u_opt_unc} and can be arbitrarily chosen.
The unconstrained value function therefore satisfies,
\begin{align}
\tilde V_k(x) = \tilde V_{k+1}(f_k(x,\kappa_k(x)))+L_k(x,\kappa_k(x)).
\end{align}

Consider the optimization problem,
\begin{subequations}\label{equ:mpc_prob_mod}
\begin{align}
\min_{\hat u} &~ \sum_{k = 0}^{N-1} \Delta \tilde V_k(x_k,\hat u_k),  \\
\text{sub.~to} &~x_{k+1} = \hat f_k(x_k,\hat u_k), \\
&~(x_k,\kappa_k(x_k)+\hat u_k) \in \mathcal C_k,~\forall k \in \mathbb Z_N, \label{equ:ineq1_mod} \\
&~x_N \in \mathcal X_N, \label{equ:ineq2_mod}
\end{align}
\end{subequations}
where $x_0$ is given and,
\begin{equation}\label{equ:deltaV}
\Delta \tilde V_k(x_k,\hat u_k) = \tilde V_{k+1}(\hat f_k(x_k,\hat u_k))+\hat L_k(x_k,\hat u_k)-\tilde V_k(x_k),
\end{equation}
with $\hat f_k(x,\hat u) = f_k(x,\kappa_k(x_k)+\hat u)$ and $\hat L_k(x,\hat u) = L_k(x,\kappa_k(x_k)+\hat u)$. 
Note that no element used in the construction of the objective function \eqref{equ:deltaV} depends on constraints.
We are now ready to state the main result.

\begin{thm}\label{thm:1}
Let the sequence of control inputs $\hat u^*$ 
solve the optimization problem \eqref{equ:mpc_prob_mod}. Then the sequence of control inputs $u^*$ 
satisfying,
\begin{equation}
u_k^* = \kappa_k(\hat x_k^*) + \hat u_k^*,
\end{equation}
for all $k \in \mathbb Z_N$, solves the optimization problem \eqref{equ:mpc_prob},
where $\hat x_{k+1}^* = \hat f_k(\hat x_k^*,\hat u_k^*)$  and $\hat x_0^* = x_0$.
\end{thm}
\begin{pf}
The sequence $\hat u^*$ 
satisfies,
\begin{subequations}
\begin{equation}\label{equ:u_mod_satisfies}
\hat u_k^* \in \underset{u \in \bar{\mathcal C}_k(\hat x_k^*)-\{\kappa_k(\hat x_k^*)\}}{\arg\min} \hat V_{k+1}(\hat f_k(\hat x_k^*,u)) + \Delta \tilde V_k(\hat x_k^*,u), 
\end{equation}
where $\hat V_k$ satisfies the Bellman equation,
\begin{equation}\label{equ:deltaV_Bellman}
\hat V_k(x) = \min_{u \in \bar{\mathcal C}_k(x)-\{\kappa_k(x)\}} \hat V_{k+1}(\hat f_k(x,u)) + \Delta \tilde V_k(x,u), 
\end{equation}
\end{subequations}
with domain $\operatorname{Proj}_{\mathbb R^n}\left(\bar{\mathcal C}_k - \{(0,\kappa_k(x))\}\right) = \operatorname{Proj}_{\mathbb R^n}\bar{\mathcal C}_k = \bar{\mathcal X}_k$, and $\hat V_N = 0|_{\mathcal X_N}$.

Let $\check V_k = \tilde V_k + \hat V_k$,~$k\in\mathbb Z_{N+1}$. Then, according to \eqref{equ:deltaV},
\begin{equation}\label{equ:long_express}
\tilde V_k(x) + \hat V_{k+1}(\hat f_k(x,u)) + \Delta\tilde V_k(x,u) = \check V_{k+1}(\hat f_k(x,u)) + \hat L_k(x,u).
\end{equation}
According to \eqref{equ:deltaV_Bellman} and \eqref{equ:long_express}, 
\begin{align}
\check V_k(x) &= \min_{u \in \mathcal C_k(x)-\{\kappa_k(x)\}} \check V_{k+1}(\hat f_k(x,u)) + \hat L_k(x,u) \nonumber \\
&= \min_{u \in \mathcal C_k(x)-\{\kappa_k(x)\}} \check V_{k+1}(f_k(x,\kappa_k(x)+u))
+ L_k(x,\kappa_k(x)+u). \label{equ:bell_check}
\end{align}
Fix $k \in \mathbb Z_N$ and assume $\check V_{k+1} = V_{k+1}$. Comparing \eqref{equ:bell1} to \eqref{equ:bell_check}, we see that they are equivalent and therefore $\check V_k = V_k$. Since $\hat V_N = 0$, then $\check V_N = 
V_N$. Therefore $\check V_k = V_k$ for all $k \in \mathbb Z_{N+1}$.
According to \eqref{equ:u_mod_satisfies} and \eqref{equ:long_express}, 
\begin{align}
\hat u_k^* &\in  \underset{u \in \bar{\mathcal C}_k(\hat x_k^*)-\{\kappa_k(\hat x_k^*)\}}{\arg\min} \check V_{k+1}(\hat f_k(\hat x_k^*,u)) + \hat L_k(\hat x_k^*,u) \nonumber \\
&= - \kappa_k(\hat x_k^*) + \underset{u \in \bar{\mathcal C}_k(\hat x_k^*)}{\arg\min}~ V_{k+1}(f_k(\hat x_k^*,u)) + L_k(\hat x_k^*,u), \label{equ:u_express}
\end{align}
in which the second expression was obtained by performing a change of variables $\kappa_k(\hat x_k^*) + u \mapsto u$.

Fix $k \in \mathbb Z_N$ and assume $\hat x_k^* = x_k^*$. Then the minimizer expressions in \eqref{equ:u_opt} and \eqref{equ:u_express} are equivalent, implying that there exists $u_k^*$ minimizing \eqref{equ:u_opt} where $u_k^* = \kappa_k(\hat x_k^*) + \hat u_k^*$ 
and $\hat x_{k+1}^* = f_k(\hat x_k^*,\kappa_k(\hat x_k^*)+\hat u_k^*) = f_k(x_k^*,u_k^*) = x_{k+1}^*$. Since $\hat x_0^* = x_0 = x_0^*$, we deduce that there exists a sequence $u^*$ 
solving \eqref{equ:mpc_prob} and satisfying $u_k^* = \kappa_k(\hat x_k^*) + \hat u_k^*$ for all $k \in \mathbb Z_N$. \qed 
\end{pf}

\section{Linear MPC optimization problem}

In practical application, the result of Theorem \ref{thm:1} is most useful in instances where there exists an analytical solution to the unconstrained version of the optimization problem \eqref{equ:mpc_prob}, or one in which the feedback law $\kappa_k(x_k)$ is conveniently parametrizable in terms of the state $x_k$. This is in particular true in the case of LQR and we consider this case in further detail by assuming a locally linear-quadratic structure to the problem \eqref{equ:mpc_prob},
\begin{equation}\label{equ:lin_approx}
\begin{split}
f_k(x,u) &= A_kx + B_ku + o(\|(x,u)\|), \\
L_k(x,u) &= \frac{1}{2}\begin{bmatrix}x^{\text T} & u^{\text T}\end{bmatrix}\begin{bmatrix}Q_k & N_k \\ N_k^{\text T} & R_k\end{bmatrix}\begin{bmatrix}x \\ u\end{bmatrix} + o(\|(x,u)\|^2), \\
V_f(x) &= \frac{1}{2}x^{\text T}P_fx + o(\|x\|^2).
\end{split}
\end{equation}
We assume $P_f$ is positive definite and introduce,
\begin{subequations}
\begin{align}
K_{k} &= -(R_{k}+B_{k}^{\text T}P_{k+1}B_{k})^{-1}(B_{k}^{\text T}P_{k+1}A_{k}+N_{k}^{\text T}), \label{equ:K} \\
P_{k} &= A_{k}^{\text T}P_{k+1}A_{k} + Q_{k} + (A_{k}^{\text T}P_{k+1}B_{k}+N_{k})K_{k}, \label{equ:P}
\end{align}
\end{subequations}
defined for $k \in \mathbb Z_N$, with $P_N = P_f$, and assume that $P_k$ and $R_{k}+B_{k}^{\text T}P_{k+1}B_{k}$ are positive definite for all $k \in \mathbb Z_N$. 


\begin{lem}
Consider the optimization problem \eqref{equ:mpc_prob} and assume that $f_k$, $L_k$, and $V_f$ are of the forms given in \eqref{equ:lin_approx}. Then,
\begin{equation}
\Delta \tilde V_k(x,\hat u) = \frac{1}{2}\hat u^{\text T}(R_k + B_k^{\mathrm T}P_{k+1}B_k)\hat u + o(\|(x,\hat u)\|^2).
\end{equation}
\end{lem}
\begin{pf}
Fix $k \in \mathbb Z_N$ and assume $\tilde V_{k+1}(x) = \frac{1}{2}x^{\text T}P_{k+1}x + o(\|x\|^2)$. Then,
\begin{equation}\label{equ:linV_min}
\tilde V_{k}(x) = \min_{u} \frac{1}{2}(Ax+Bu)^{\text T}P_{k+1}(Ax+Bu) + \frac{1}{2}\begin{bmatrix}x^{\text T} & u^{\text T}\end{bmatrix}\begin{bmatrix}Q_k & N_k \\ N_k^{\text T} & R_k\end{bmatrix}\begin{bmatrix}x \\ u\end{bmatrix} + o(\|(x,u)\|^2).
\end{equation}
As shown by \cite{dorato71}, the unique minimizer of \eqref{equ:linV_min} is given by $u = K_kx + o(\|x\|)$ and,
\begin{equation}\label{equ:linV}
\tilde V_{k}(x) = \frac{1}{2}x^{\text T}P_{k}x + o(\|x\|^2).
\end{equation}
Since $V_N = V_f$, we can deduce that \eqref{equ:linV} is true for all $k \in \mathbb Z_N$. According to \eqref{equ:deltaV},
\begin{equation*}
\Delta \tilde V_k(x,\hat u) = \frac{1}{2}\begin{bmatrix}x^{\text T} & \hat u^{\text T}\end{bmatrix}
\begin{bmatrix}
\hat Q_{k} & \hat N_{k} \\ \hat N_{k}^{\text T} & \hat R_{k}
\end{bmatrix}\begin{bmatrix}x \\ \hat u\end{bmatrix} + o(\|(x,\hat u)\|^2),
\end{equation*}
where,
\begin{align*}
\hat Q_k &= \hat A_k^{\text T}P_{k+1}\hat A_k + Q_k + K_k^{\text T}N_k^{\text T} + N_kK_k + K_k^{\text T}R_kK_k - P_k, \\
\hat N_k &= \hat A_k^{\text T}P_{k+1} B_k + N_k + K_k^{\text T}R_k, \\
\hat R_k &= B_k^{\text T}P_{k+1} B_k + R_k,
\end{align*}
and $\hat A_k = A_k+B_kK_k$.
Note that \eqref{equ:K} and \eqref{equ:P} imply that $\hat N_k = 0$ and $\hat Q_k = 0$, respectively.
\qed
\end{pf}

Consider the following optimization problem,
\begin{subequations}\label{equ:mpc_prob_lin}
\begin{align}
\min_{\hat u} &~ \frac{1}{2}\sum_{k = 0}^{N-1} \hat u_k^{\text T}(R_k+B_k^{\text T}P_{k+1}B_k)\hat u_k,  \label{equ:cost_lin}\\
\text{sub.~to} &~x_{k+1} = \hat A_kx_k+B_k\hat u_k, \\
&~(x_k,K_kx_k+\hat u_k) \in \mathcal C_k,~\forall k \in \mathbb Z_N, \label{equ:ineq1_lin} \\
&~x_N \in \mathcal X_N, \label{equ:ineq2_lin}
\end{align}
\end{subequations}
with $x_0$ given. 
The following result is a straightforward application of Theorem \ref{thm:1}.
\begin{cor}\label{cor:3}
Consider \eqref{equ:mpc_prob} and assume that $f_k$, $L_k$, and $V_f$ are of the forms given in \eqref{equ:lin_approx} with no residual terms, i.e., $o = 0$.
Let the sequence of control inputs $\hat u'$ 
solve \eqref{equ:mpc_prob_lin}. Then the sequence of control inputs $u'$ 
satisfying,
\begin{equation}
u_k' = K_k\hat x_k'+\hat u_k',
\end{equation}
for all $k\in\mathbb Z_N$, solves the optimization problem \eqref{equ:mpc_prob}, where $\hat x_{k+1}' = \hat A_k\hat x_k'+B_k\hat u_k'$ and $\hat x_0' = x_0$. \qed
\end{cor}


\section{Constraint-separation principle}

We refer to the result of Corollary \ref{cor:3} as a constraint-separation principle, because it separates constraint enforcement from stabilization in MPC. 
The result implies that any locally linear-quadratic, open-loop MPC problem can be locally restructured as a closed-loop, constraint-enforcing, minimum-norm projection problem, 
where the feedback gain of the closed-loop controller is the optimal gain obtained by solving the unconstrained, open-loop MPC problem.

Importantly, this equivalence shows that linear MPC can be interpreted as an application of a constraint enforcement scheme to an already precompensated system, an interpretation that is associated with reference and command governors \citep{rg_survey} and, in addition to being applied to reference commands, is what supposedly, significantly distinguishes these methods from MPC. It also shows that the extended command governor (ECG) \citep{gilbertong11}, if applied as an offset to a control input similarly to the approach taken by \citet{rossiter10}, can be viewed as a generalization of linear MPC, because the problem \eqref{equ:mpc_prob_lin} is a special case of the conventional ECG optimization with reference input set to $0$.\footnote{Details establishing the link between ECGs and MPC are available in Section 5 of \citet{gilbertong11}.} 

The result also provides a method by which to choose the MPC penalty function when applying the decomposition technique to precompensated systems. For example, given a sequence of gains $K_k$, the penalty matrix should be set to $R_k+B_k^{\text T}P_{k+1}B_k$, where $R_k$ and $P_k$ 
satisfy the solution of an inverse LQR problem obtained, for example, using methods of \citet{kreindler72,dicairano10}. In the nonlinear case, a solution to a general, inverse optimal control problem could be used to determine a cost function for an MPC controller with the cost function $\Delta V_k$ being determined according to the results of the main theorem.



Additionally, as explored by \citet{kouvaritakis2000,kouvaritakis02,rossiter_book}, the result provides a potential simplification to numerical approaches to MPC problems. As discussed by \cite{numerical_mpc}, most conventional numerical approaches to nonlinear MPC problems either apply a sequential quadratic programming approach, where the corollary is obviously useful, or an interior point approach, where the main theorem is more useful since precompensation at least allows for better numerical conditioning as it prevents the closed-loop maps $\hat f_k$ from blowing up.
To the best of our knowledge, most solvers applied to MPC, although allowing for the option to precompensate according to the optimal feedback, e.g., \citet{mpt}, do not fully take advantage of the separation principle presented. For example, they do not decompose the local problem into what it essentially consists of: an LQR problem and a minimum-norm projection onto a closed set. They instead allow for a user to precompensate according to any feedback when the choice of optimal feedback has greater utility, as it guarantees no error when constraints are inactive. 

\begin{figure}[tbp]
\centering
\begin{tikzpicture}[auto, node distance=2cm,>=latex']
    \node [input, name=input] {};
    \node [block, right of=input, text width=0.7in,align=center] (controller) {optimal feedback controller};
    \node [sum, right of=controller] (sum) {};
    \node [block, right of=sum, text width=0.7in,align=center] (system) {open-loop plant};
    \node [block, above of=controller, node distance=1.4cm, text width=0.7in, align=center] (mpc) {closed-loop MPC};
    \draw [->] (controller) -- node[name=Kx] {$\kappa(\hat x)$} (sum);
    \draw [->] (sum) -- node[name=u,pos=0.4] {$u$} (system);
    \draw [->] (mpc) -| node[name=uhat,pos=0.25] {$\hat u$} (sum);
    \node [output, right of=system] (output) {};
    \node [block, below of=sum] (measurements) {observer};
    
    \draw [->] (input) -- node [pos=0.45] {$\hat x$} (controller);
    \draw [->] (system) -- node [name=y,pos=0.4] {$y$}(output);
    \draw [->] (y) |- (measurements);
    \draw [-] (measurements) -| 
        node [near end] {} (input);
    \draw [->] (input) |- (mpc);
\end{tikzpicture}
\caption{MPC decomposed into an optimal feedback controller and a constraint-enforcing scheme}\label{fig:MPC_cl}
\end{figure}
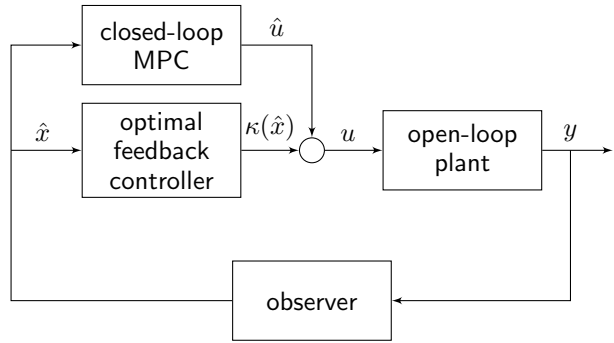

The discussion thus far has focused on the special case where the MPC problem is locally linear-quadratic. We note that the main result is applicable more generally, as represented by the schematic in Fig.~\ref{fig:MPC_cl}, which we have shown to be equivalent to Fig.~\ref{fig:MPC_ol}. The constraint-separation principle does not necessarily hold in the general case, 
when the problem is not locally linear-quadratic, e.g., when there does not exist a continuously differentiable stabilizing feedback $\kappa_k$ or the cost function $V_k$ is not continuously differentiable; in this case, it is not guaranteed that 
the penalty function $\Delta \tilde V_k$ is locally independent of $x$. Nevertheless, the generality of the main result is remarkable: The result can be applied as long as a solution is known to the unconstrained problem, with a few minor, technical conditions corresponding to constraints and penalty function. However, note that constraint-enforcement does not necessarily simplify to a minimum-norm projection in this case.

Now note the practical use of decomposition in certifying MPC controllers. A decomposed controller is more straightforward to certify because an unconstrained problem is simpler than the same problem with constraints; it is therefore easier to certify the stabilizing component of the decomposed MPC problem. This makes it easier to ensure stability and attractiveness during online operation, limiting the need for a complex certification process, such as that of \citep{garoche18}, to the constraint-enforcement component.


Taken together, the above discussion implies that there is a strong desire to often, if not always, decompose design in the manner derived in this manuscript, whenever one is able to effectively parametrize the feedback controller $\kappa_k$ and the closed-loop state update $\hat f_k$. This approach can, for example, improve design of neural-net-based MPC \citep{saintdonat91} by decomposing the problem into the solution of a simpler, unconstrained, approximate dynamic programming problem \citep{bertsekas_dp}, and a more difficult, constrained optimization problem. 

As a matter of course, this discussion has only been able to superficially consider the practical use of decomposition in nonlinear MPC. We feel, however, that it represents a promising direction for future research.





\section{Conclusion}
In this brief, we derived a constraint-separation principle for MPC problems. The results show that MPC problems can be decomposed into the solution of an unconstrained, open-loop problem and a constrained, closed-loop problem without a terminal cost, which may simplify MPC problems when the unconstrained, stabilizing feedback can be represented explicitly by a closed-form solution or a parametrization. This is particularly true for linear MPC problems, for which the stabilizing feedback is given as the solution to the well-known LQR problem. 
It is significant because it shows the equivalence of designing MPC in a two-step approach, which first stabilizes the system and then implements constraint protection in the outer-loop. It is also significant because it can be used to improve numerical solutions to both linear and nonlinear MPC problems. 

\bibliography{autosam}           



\end{document}